\documentclass[11pt,a4paper]{article}
\usepackage{amsmath,amsthm,amssymb}
\usepackage{graphicx}
\usepackage{tikz}
\usepackage{pgfplots}
\pgfplotsset{compat=1.13} 
\usepackage{verbatim}
\usepackage{color}
\usepackage{subfigure}
\usepackage{float}
\usepackage{authblk}
\usepackage[colorlinks]{hyperref}

\usepackage{times}

\usepackage{textcomp}
\usepackage{eurosym}
\theoremstyle{plain}
\newtheorem{thm}{Theorem}[]
\newtheorem{lem}[thm]{Lemma}
\newtheorem{cor}[thm]{Corollary}

\theoremstyle{definition}

\newtheorem{example}{Example}[]

\theoremstyle{remark}

\renewcommand{\Re}{\operatorname{Re}}

\begin{document}

\title{Exact eigenvalue assignment of linear scalar systems with single delay using Lambert W function}
\author[1]{Huang-Nan Huang}
\author[2]{Chew Chun Yong}

\affil[1]{Department of Applied Mathematics, Tunghai University, Taichung, Taiwan (E-mail: nhuang@thu.edu.tw)}
\affil[2]{Department of Mathematical \& Actuarial Sciences, University Tunku Abdul Rahman, Bandar Sungai Long, Selangor, Malaysia (E-mail: chewcy@utar.edu.my)}

\date{March 15, 2018}

\maketitle
\begin{abstract}
Eigenvalue assignment problem of a linear scalar system with a single discrete delay is analytically and exactly solved. The existence condition of the desired eigenvalue is established when the current and delay states are present in the feedback loop. Design of the feedback controller is then followed. Furthermore, eigenvalue assignment for the input-delay system is also obtained as well. Numerical examples illustrate the procedure of assigning the desired eigenvalue.

\vskip1em \noindent \textbf{2010 AMS Classification:}
39B82, 93B55
\end{abstract}

\section{Introduction}
During recent decades, the stabilization and control of linear systems with delays are extensively studied, for example, the spectrum (eigenvalues) assignment for linear delay systems in 1978 by Olbrot \cite{Olbrot}. Recently, an approach for the solution of linear time-delay system is based on the Lambert W function proposed by Asl and Ulsoy \cite{Asl}. Hence the robust stability as well as related topics to design the feedback controller are well established \cite{Yi2010a} and reference therein. A good introduction about the Lambert W function is given by Coreless \emph{et. al.} \cite{Coreless} and this function possesses many applications within these two decades.

Eigenvalue assignment for delay systems with single delay via Lambert W function is first developed by Yi \emph{et. al.} \cite{Yi2010b} to assign the rightmost eigenvalue of the delay system to a predefined (desired) location in order to stabilize the system, but unfortunately only a real or real part of the rightmost eigenvalue can be assigned for the scalar case. Alternatively, Shinozaki \cite{Shinozaki07} discussed the assignment of the complex eigenvalue to the largest one of a scalar single delay system with the output of a complex feedback gain which is not realistic. Later an analytic eigenvalue assignment method is proposed for scalar and some special delay systems \cite{Nemati}. All these studies design the controller by feedback only the current state and no condition is drawn on the value of desire eigenvalue such that the feedback controller always exists. On the other hand, although a more general time-delay system can be analyzed by using matrix Lambert W function, but the rightmost eigenvalue can be computed not by using the principal branch which contradicts the main proposition of this method \cite{Cepeda}.
 
Avoiding using the matrix Lambert W function, only a scalar system with a single delay is considered in this paper. It mainly focuses on deriving the existence condition of the feedback controller related to assign the rightmost eigenvalue of the system to the desired value. The formula to compute feedback gains for the current and delay states is then obtained. Furthermore, an input-delay type system is also considered. Two examples are provided to demonstrate the associated idea. The result consolidate further studies on the eigenvalue assignment of linear systems with multiple delays via the Lambert W function approach.

\section{Lambert W functions}
\label{sec2}

\begin{figure}[htbp]
    \centering
\begin{tabular}{cc}
    \includegraphics[width=6.6cm]{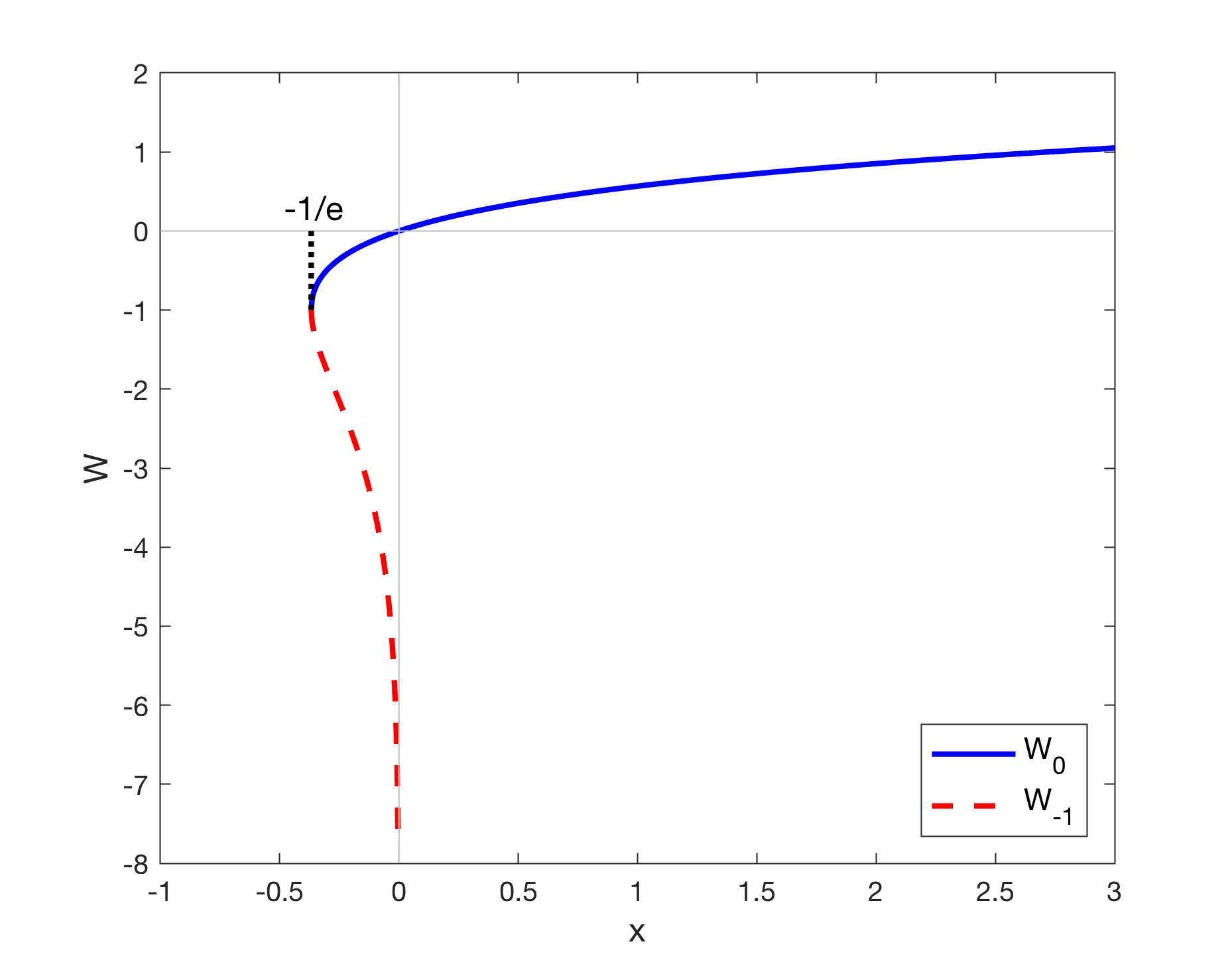}
    &
    \includegraphics[width=7cm]{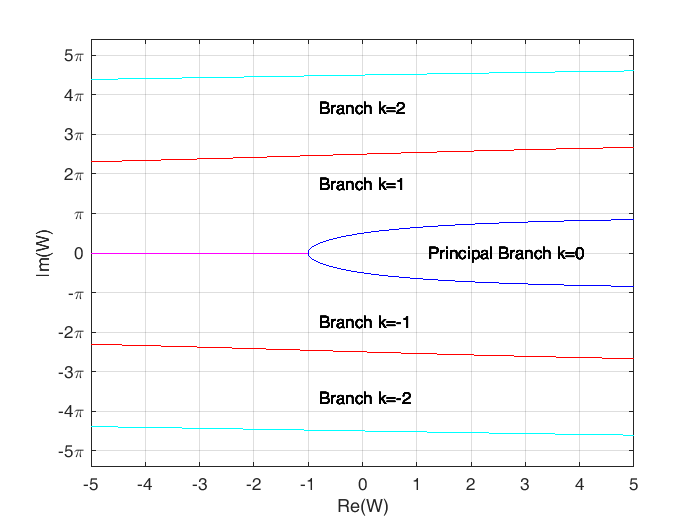}\\
    (a) $W_0(x)$,$W_{-1}(x)$,$x\in\mathbb{R}$. & (b) Ranges of $W_{k}$, $k=0,\pm1,\pm2$.
\\
\end{tabular}
    \caption{The range of Lambert W function.}
    \label{fig:W_range}
\end{figure}

The Lambert W function is defined as a complex multivalued function which has infinite number of branches, $W_k(x)$, where $k=0, \pm 1,  
\pm 2,\cdots,\pm\infty$,
(regard $W_\infty$ and $W_{-\infty}$ as fixed mappings), such that
\[
W_k(z) e^{W_k(z)} =z,\quad z\in\mathbb{C}.
\]
For any $x\in\mathbb{R}$, when $-\frac{1}{e}\le z=x < 0$, the \emph{principal branch} $W_0(x)$ satisfies $W_0(x)\geq{-1}$ and the $-1$ branch  $W_{-1}(x)$ satisfies $W_{-1}(x)\leq 
-1 $. 
By partitioning the $z$-plane with horizontal boundaries $z=(2k+1)\pi+i~0 $  for $k\in\mathbb{Z}$ (where $i=\sqrt{-1}$), the ranges of branches of $W_k(z)$ are images of the $z$ between branch cuts in the $z$-plane. $W_0$ has a branch cut linking to $W_1$ and $W_{-1}$ which is defined as $BC=\{x+i~0| x\in(-\infty,-1/e)\}$, i.e., $W_0(BC)$ is the boundary between ranges of $W_0$ and $W_1$ and so as $W_{-1}(BC)$ for $W_{-1}$ and $W_0$. The range of branches as well as its real counterpart of this function are shown in Figure \ref{fig:W_range}. 

An important property of the Lambert W function is given by \cite{Shinozaki07,Huang}
\begin{lem}
	\label{lem:1}
	Let $z\in\mathbb{C}$. Then  $\overline{W_k(z)} = W_{-k}(\bar{z})$
		and $ \max\limits_{k=0, \pm1,\pm2,\cdots,\pm\infty} \Re[ W_{k}(z)]=\Re[ W_{0}(z)]$.
\end{lem}

\section{Main Result}

We here consider a scalar delay system with an exogenous input from environment: 
\begin{equation}
\begin{split}
& \dot{x} (t) =a x(t)+ a_{1d} x(t-h)+b u(t),\quad h>0,\\
& x(0)  = x_0,\quad    x( \tau ) =\phi(\tau),\quad -h\le \tau<0,
\end{split}  
\label{eqn:1}   
\end{equation}   
where $a, a_{1d}, b\neq 0, x_0,~h\in\mathbb{R}$, and $\phi$ is the initial function
to specify  $x(\tau)$, $\tau \in[-h,0)$.
Suppose a proportional control is proposed to stabilize the system with feedback of current and delay states:
\begin{equation}
   u(t)=k x(t) +k_{1d} x(t-h)
   \label{eqn:controller}
\end{equation}
with two parameters $k,~ k_{1d} \in \mathbb{R}$ to be designed.
The closed-loop system is then described by
\begin{equation}
\dot{x} =(a+bk) x(t)+ (a_{1d}+b k_{1d}) x(t-h)  \triangleq \alpha x(t) + \beta x(t-h)
\label{eqn:single-close}
\end{equation}
where $\alpha=a+b k$ and $\beta=a_{1d}+b k_{1d}$ which are both real numbers. The characteristic equation of the closed-loop given by \cite{Huang}:
\begin{equation}
    s-\alpha=\beta e^{-s h},
    \mbox{ or equivalently, }
    s-(a+bk) = (a_{1d}+b k_{1d}) e^{-s h}
    \label{eqn:quasi-poly}
\end{equation}
whose roots (also known as \emph{eigenvalue} of the system) are
expressed by using Lambert W function
\begin{equation}
   s_k=\alpha+\frac{1}{h}W_k(\beta h e^{-\alpha h}),
   \quad k=0, \pm 1, \pm 2,\ldots,\pm \infty.
   \label{eqn:roots}
\end{equation}
From Lemma \ref{lem:1}, the system
(\ref{eqn:single-close})
is stable iff the real part of the rightmost eigenvalue, $s_0$, is negative.

The control design about a system with only input delay Whose equation is given by
\begin{equation}
   \dot{x} =a x(t) + b u(t-h),\quad x(0)  = x_0, \mbox{ and } h>0,
   \label{eqn:input-delay}   
\end{equation}
is much more a demanding challenge since the corresponding input-output operator is not a compact one. 
Applying the state feedback controller
$u(t) =k x(t)$, then the associated closed-loop system becomes
$
\dot{x} = a x(t)+b k x(t-h)
$
which is also of the form (\ref{eqn:single-close}) with $\alpha=a$ and $\beta = b k$.
Thus, we only need to focus on the system (\ref{eqn:single-close}).

How to assign the rightmost  eigenvalue $s_0$ of the closed-loop system to a desired location?  
Suppose the desired location is denoted by $S_{0,des}\in\mathbb{C}$ (with positive imaginary part) but not a real number, we want to calculate real parameters $k$ and $k_{1d}$ (i.e., adjust the values of $\alpha$ and $\beta$) such that $s_0=S_{0,des}$.  As discussed in \cite{Huang},
$\bar{S}_{0,des}$ must be the eigenvalue $s_{-1}$ of the system belonged to the range of $W_{-1}$.
Let $S_{0,des}=u+i v$ and substituting $s$ in (\ref{eqn:quasi-poly}) leads to
$u+i v-\alpha = \beta e^{-(u+i v)h}$, or equivalently,
\begin{equation}
    \begin{cases}
        u-\alpha  = \beta e^{u h}\cos v h, \\
        v =-\beta e^{u h} \sin v h.
    \end{cases}
    \label{eqn:exist1}
\end{equation}
Combine to yield $(u-\alpha)h=-v h \cot v h$, i.e., $(S_{0,des}-\alpha)h = (u-\alpha) h+ i v h \in W_0(BC)$, i.e., $S_{0,des}-\alpha$ lies in the boundary between the ranges of $W_0$ and $W_1$. From (\ref{eqn:exist1}) we have $\beta =- v e^{-u h} \csc v h$ (note $v\neq 0$) and $\alpha=u+v \cot v h$ and then the associated control feedback gains are given by
\begin{equation}
    \begin{cases}
       \displaystyle  k =\frac{u+v \cot v h -a}{b}, \\
      \displaystyle   k_{1d} = -\frac{v e^{-u h} \csc v h+a_{1d}}{b}.
    \end{cases}
    \label{eqn:gain1}
\end{equation}

What happens when one of the states in (\ref{eqn:controller}) is not included in the feedback? Suppose the current state is not feedback ($k=0$), the rightmost eigenvalue is assignable if $(S_{0,des}-a) h\in W_0(BC)$, or equivalently, $a=u+v \cot v h$, and the delay state feedback gain is then given by
\begin{equation}
    k_{1d}=\frac{ (S_{0,des}-a) e^{S_{0,des}h}-a_{1d}}{b}.
    \label{eqn:gain2}
\end{equation}
Furthermore, suppose the delay state is not used ($k_{1d}=0$), we obtain $S_{0,des}-a_{1d}e^{-S_{0,des}h} =  a+bk\in\mathbb{R}$ from (\ref{eqn:quasi-poly}), i.e., $S_{0,des}$ must satisfy the condition $a_{1d}+v e^{u h} \csc vh=0$. The current state feedback gain is then described by
\begin{equation}
    k=\frac{ S_{0,des}-a -a_{1d} e^{-S_{0,des}h}}{b}.
    \label{eqn:gain3}
\end{equation}

On the other hand, if $S_{0,des}\in\mathbb{R}$, then (\ref{eqn:roots}) becomes
$
S_{0,des} = \alpha + \frac{1}{h} W_0(\beta h e^{-\alpha h})
$, i.e.,
$S_{0,des}\ge \alpha -\frac1{h}$ while
$\beta h e^{-\alpha h}\ge -\frac1e$, due to the range of $W_{0}$, lies in the right-hand side of the vertical line $u=-1$. Unless the system is an input-delay system, we always can find an $\alpha$ (or $k$) such that $\alpha \le S_{0,des}+\frac{1}{h}$ and $\beta=(S_{0,des}-\alpha) e^{S_{0,des}h}$ with corresponding $k=(\alpha-a)/b$ and $k_{1d}=(\beta-a_{1d})/b$, i.e.,
\begin{equation}
    \begin{cases}
    \displaystyle k\le \frac{S_{0,des}-a+1/h}{b}, \\
    \displaystyle k_{1d} = \frac{[S_{0,des}-(a+bk)]e^{S_{0,des}h}-a_{1d}}{b}
    \ge \frac{e^{S_{0,des}h}}{hb}.
    \end{cases}
    \label{eqn:gain4}
\end{equation}
Suppose the current state is not present in the feedback ($k=0$), the rightmost eigenvalue is assignable if $S_{0,des}\ge a-\frac1h$, and the delay state feedback gain is still given by (\ref{eqn:gain2}).
Similarly, when the delay state is not used ($k_{1d}=0$), we have $S_{0,des}-a_{1d}e^{-S_{0,des}h} =a +b k$ which is always achievable by the state feedback gain from (\ref{eqn:gain3}).

An input delay system (\ref{eqn:input-delay})is assignable to any complex $S_{0,des}$ if
$(S_{0,des}-a)h\in W_0(BC)$, i.e., $(S_{0,des}-a)h e^{(S_{0,des}-a)h} =z$ for some real number $z<-\frac1e$. Then the associated real feedback gain for the controller $u=k x(t)$ is determined by
\begin{equation}
  k=\frac{S_{0,des}-a}{b}e^{ S_{0,des}h}.
\label{eqn:gain-input-delay}
\end{equation}
Or when $S_{0,des}$ is real, it must satisfy $S_{0,des}\ge a-\frac1h$ and the feedback gains is still given by (\ref{eqn:gain-input-delay}) which is the same as the result presented in \cite{Nemati}.

From the above derivation, the following result is asserted:
\begin{thm}\label{thm:2}
Suppose the system \emph{(\ref{eqn:1})} in not an input-delay system, the following statements hold:
\begin{enumerate}
    \item[(i)] For a given $S_{0,des}=u+i v\in\mathbb{C}$, the rightmost eigenvalue of the closed-loop system \emph{(\ref{eqn:single-close})} can be assignable to any desired location $S_{0,des}$ via the controller \emph{(\ref{eqn:controller})} with both current- and delay-state feedback gains defined by \emph{(\ref{eqn:gain1})}.\\
    Furthermore, if the current or delay state is not included in the feedback loop, the $S_{0,des}$ must satisfy the condition $a = u+v \cot v h$ or $a_{1d}=-v e^{u h} \csc vh$ such that the associated gain is described by \emph{(\ref{eqn:gain2})} or \emph{(\ref{eqn:gain3})}, respectively.
   
    \item[(ii)]
    For a given $S_{0,des}\in\mathbb{R}$, the rightmost eigenvalue of the system \emph{(\ref{eqn:single-close})} can be assignable to any desired location $S_{0,des}$ via the controller \emph{(\ref{eqn:controller})} with feedback gains defined by \emph{(\ref{eqn:gain4})}. Furthermore, if the current or delay state is not included in the feedback loop, the $S_{0,des}$ must satisfy the condition $S_{0,des}\ge a-\frac1h$ or no constraint such that the associated gain is still described by \emph{(\ref{eqn:gain2})} or \emph{(\ref{eqn:gain3})}, respectively. 
\end{enumerate}
\end{thm}
\begin{cor}\label{cor:3} For an input-delay system \emph{(\ref{eqn:input-delay})}, if $S_{0,des}-\alpha$ belongs to the upper boundary on the range of $W_0$ or $[-1,-\infty)$, a real feedback gain $k$ through \emph{(\ref{eqn:gain-input-delay})} is obtained. 
\end{cor}

\begin{example} \label{exmp:1}
Consider the system (\ref{eqn:1}) with the given data set $a=1$, $a_{1d}=-1$, $b=1$, and $h=1$. Determine the values of $k$ and $k_{1d}$ such that the rightmost closed-loop eigenvalue $s_0$ is located at $S_{0,des}=−0.092484+1.9973i$, $−0.60502+1.7882i$, or $−1$, respectively.

By Theorem \ref{thm:2} this system is eigenvalue-assignable to three desired eigenvalues with the feedback gains given by (\ref{eqn:gain1}) for the first two or (\ref{eqn:gain4}) for the third one. Table \ref{tab:1} shows the computational result for these parameters. And the corresponding variation of characteristic roots of the closed-loop system before and after the eigenvalue assignment are also shown in Table \ref{tab:2} by using the method proposed in \cite{Huang}. For the third eigenvalue, the closed-loop system becomes a non-delay system which has only one eigenvalue $-1$ as expected.

\begin{table}[htb]
   \centering   
   \caption{State feedback gains of the controller with respect to three different eigenvalue assignments.}
   \label{tab:1}
   \begin{tabular}[t]{|c|c|c|c|}
    \hline 
    $S_{0,des}$ & $-0.092484+1.99730i$ & $-0.60502+1.78820i$ &  $-1.0+0i$ \\ 
    \hline
    $k$ & $-2$ & $-2$ & $-2$ \\ 
    \hline
    $k_{1d}$ & $-1$ & $~~0$ & $~~1$ \\
    \hline
    \end{tabular}
\end{table}

\begin{table}[htb]
    \centering
    \caption{The variation of characteristic roots before and after eigenvalue assignment.}
    \label{tab:2}
    \begin{tabular}[t]{|c||c|c|}
    \hline 
    $a=1$, $a_{1d}=-1$ & $k=-2$, $k_{1d}=-1$ & $k=-2$, $k_{1d}=0$\\ 
    \hline
    \hline
    $s_{3,-3}=-3.02630\pm 20.2238i$ & $s_{3,-4}=-2.32231\pm 20.3555i$ & $s_{3,-4}=-3.01658\pm 20.3214i$ \\ 
    \hline 
    $s_{2,-2}=-2.66407\pm 13.8791i$ & $s_{2,-3}=-1.95315\pm 14.0695i$ & $s_{2,-3}=-2.64736\pm 14.0202i$  \\
    \hline 
     $s_{1,-1}=-2.08880\pm 7.46150i$ & $s_{1,-2}=-1.36300\pm 7.80750i$& $s_{1,-2}=-2.05280\pm 7.71840i$ \\
    \hline
    $s_0=0$ & $s_{0,-1}=-0.092484\pm 1.99730i$ & $s_{0,-1}=-0.60502 \pm 1.78820i$ \\ 
    \hline 
    \end{tabular}
\end{table}
        
\end{example} 

\vspace{-2em}
\section{Conclusion}
The sufficient condition on the solving the eigenvalue assignment of a linear scalar systems with single delay is analytically developed. The feedback controller is designed by using both the current and delay states. Furthermore, similar sufficient condition for input-delay systems is also depicted. The gains to associate feedback states are computed accordingly such that the closed-loop system behaves the same as expected. Examples are presented for illustrate purpose. This result provides more flexibility in real-world applications than previous studies.

\subsection*{Acknowledgement}
This work was partially supported under the grant No. MOST 106-
2115-M-029-004.

\end{document}